\documentstyle[12pt, leqno]{article}
\newtheorem{theorem}{Theorem}[section]

\newtheorem{proposition}[theorem]{Proposition}
\newtheorem{corollary}[theorem]{Corollary}
\newtheorem{lemma}[theorem]{Lemma}

\def \proof {\noindent {\bf Proof.}\ \ }
\def \endproof {{\mbox{}\nolinebreak\hfill\rule{2mm}{2mm}\par\medbreak} }

\def \ll {\langle}
\def \rr {\rangle}

\def \R {{\bf R}}
\def \a {\alpha}
\def \b {\beta}
\def \e {\varepsilon}
\def \d {\delta}
\def \k {\kappa}
\def \l {\lambda}
\def \s {\sigma}
\def \t {\tau}
\def \w {\omega}
\def \E {{\bf E}}
             
\def \R {{\bf R}}
\def \< {\langle}
\def \> {\rangle}
\def \ra {\rightarrow}

\def \dist {{\rm dist}}
\def \span {{\rm span}}
\def \rank {{\rm rank}}
\def \trace {{\rm trace}}

\def \HS {{\rm HS}}
\def \DR {{\rm DR}}
\def \id {{\it id}}

\def \vr {{\rm vr}}
\def \vol {{\rm Vol}}
\def \dens {\overline{\rm dens}\;}

\begin{document}
\title {John decompositions: selecting a large part}
\author {R. Vershynin\footnote{Department of Mathematics, University of Missouri, 
                               Columbia, MO 65211, USA (vershynin@math.missouri.edu)} }
\date{September 15, 1999}
\maketitle

\begin{abstract}
  We extend the invertibility principle of J.~Bourgain and L.~Tzafriri
  to operators acting on arbitrary decompositions
  $\id  =  \sum x_j \otimes x_j$, rather than on the coordinate one.
  The John's decomposition brings this result to the local theory of Banach spaces.
  As a consequence, we get a new lemma of Dvoretzky-Rogers type, 
  where the contact points of the unit ball with its maximal volume ellipsoid
  play a crucial role. This is applied to embeddings of $l_\infty^k$
  into finite dimensional spaces.
\end{abstract}


\section{Introduction}

The aim of this paper is to find a part of a John decomposition
on which a given nontrivial operator is invertible in a certain sense, 
and to apply this to the study of contact points of convex bodies. 

John decomposition of identity is by now a classical tool in the local 
theory of Banach spaces.
Suppose $X  =  (\R^n, \| \cdot \|)$ is a Banach space whose 
ellipsoid of maximal volume contained in $B(X)$ 
coincides with the unit euclidean ball. 
The John decomposition of the identity operator on $X$ is
\begin{equation}                                                             \label{Johndecomp}
  \id  =  \sum_{j = 1}^m  x_j \otimes x_j,
\end{equation}
where $x_j / \|x_j\|_X$ are some contact points of the surfaces
of $B(X)$ and the unit euclidean ball.
This celebrated theorem of F.~John has been used extensively 
over past 30 years. Recently it was interpreted as an isotropic
condition \cite{G-M}, and was generalized to a non-convex case
in \cite{G-P-T}.

It is important to know good parts of John decomposition, 
as it provides a good set of contact points of $X$, 
which can be useful for understanding the geometrical structure of $X$,
see \cite{R2}.
In the present paper we find a part of a decomposition (\ref{Johndecomp})
which preserves the orthogonal structure under action of a given 
linear operator $T$.
More precisely, if $\|T\|_{2 \ra 2}  =  1$
then there exists a subset of indices $\s$
of cardinality $|\s|  \ge  (1 - \e) \|T\|_\HS^2$ such that
the system $(T x_j)_{j \in \s}$ is $C(\e)$-equivalent
to an orthogonal basis in Hilbert space.
This result is non-trivial even for the coordinate decomposition
$\id  =  \sum e_j \otimes e_j$.
In this case it generalizes the principle of restricted invertibility
proved by J.~Bourgain and L.~Tzafriri \cite{B-Tz}. 
They considered only operators $T$ for which all norms $\|T e_j\|_2$ 
are well bounded below and proved the principle with some fixed $0 < \e < 1$.

$T$ being an orthogonal projection, we derive a new lemma of 
Dvoretzky-Rogers type. 
Suppose $P$ is an orthogonal projection in $X$ with $\rank P = k$.
Then for any $\k  <  k$
there are contact points $x_1, \ldots, x_\k$
such that setting $z_j  =  P x_j / \|P x_j\|_2$
we have

\begin{itemize}
\item{the system $(z_j)$ is $C(\k / k)$-equivalent in $l_2$-norm 
    to the canonical basis of $l_2^\k$;}
\item{$\|z_j\|_X  \ge  c \sqrt{\frac{k - \k}{n}}$  for all $j$.}
\end{itemize}

\noindent To put the result in other words, 
the orthonormal system in $Z$ guaranteed by the classical 
Dvoretzky-Rogers Lemma is essentially the normalized projections 
of contact points of $X$.
Moreover, the result holds for selfadjoint operators as well 
as for projections, the Hilbert-Schmidt norm substituting $\rank P$.
For a general operator $T$ the best lower bound for $\|T x_j\|_X$ 
is equivalent to $\frac{1}{n} |\trace T|$.

$T$ being the identity operator, we obtain a set of 
$k  >  (1 - \e) n$ contact points of $X$
which is $C(\e)$-equivalent in $l_2$-norm to the canonical 
basis in $l_2^k$.
This settles an isomorphic version of a problem of 
N.~Tomczak-Jaegermann (\cite{T-J}, p.127),
and confirms the feeling that contact points are always 
distributed fairly uniformly on the surface of the maximal volume ellipsoid
(see \cite{B1}).
Besides, this yields the proportional Dvoretzky-Rogers factorization
(with constant $C(\e)  =  \e^{c \log \e}$, 
which is however not the best known estimate).

The Dvoretzky-Rogers Lemma is proved useful in study 
of subspaces of $X$ well isomorphic to $l_\infty^k$.
The use of the refined Dvoretzky-Rogers lemma above
improves a "Gaussian" version of Alon-Milman-Talagrand Theorem 
about $l_\infty^k$-subspaces of $X$.
Let $P$ be an orthogonal projection in $X$ with $\rank P = k$.
Then there exists a subspace $Z  \subset X$ 
which is $M$-isomorphic to $l_\infty^m$ 
for $m  \ge  c k / \sqrt{n}$
and $M  =  c \sqrt{\frac{n}{k}} \ell(P)$.
The subspace $Z$ is canonically spanned by the projections
of $m$ contact points $x_j$.
Moreover, the norm on $Z$ is $M$-equivalent to 
$|||z|||  =  \max_{j \le m} | \< z , x_j \> |$.
This improves the estimates obtained by 
M.~Rudelson in \cite{R1}, and also provides information
about the position of $Z$ in $X$.
Besides, this yields a refinement of M.~Rudelson's 
result about $l_\infty^k$-subspaces in spaces
with large volume ratio.
If $\vr(X)  \ge  a \sqrt{n}$ then $X$ has a subspace $Z$ 
of dimension $m  \ge  C_1(a) \sqrt{n}$
which is $C_2(a) \log{n}$-isomorphic to $l_\infty^m$.

The extraction results about John decompositions
can be reformulated in the language of frames. 
Suppose we are given a tight frame $(x_j)$ in Hilbert space $H$,
and a norm-one linear operator $T : H  \ra  H$. 
Then there is a subsequence 
$(T x_j)_{j \in \s}$ with $|\s|  \ge  (1 - \e) \|T\|_\HS^2$ 
which is $C(\e)$-equivalent to an orthogonal basis in Hilbert space.
This theorem again can be interpreted as an extension of
the invertibility principle. 
It also generalizes results of P.~Casazza \cite{C2} and 
the author \cite{V}, who worked with the identity operator 
$T = \id$.

The rest of the paper is organized as follows. 
In \S \ref{SecPrelim} we recall some basic tools used later. 
The extraction result about John decompositions, as
well as some modifications, is proved in \S \ref{SecMain}.
Its relation to the principle of restricted invertibility
and infinite-dimensional analogs are discussed in \S \ref{SecInvert}.
Dvoretzky-Rogers type lemmas are derived from these results
in \S \ref{SecDR}.
They help to understand structure of the set of contact points. 
Applications to $l_\infty^k$-subspaces of a finite dimensional 
space are given in \S \ref{SecCube}.
Finally, in \S \ref{SecFrames} we discuss a relation of 
these results to the theory of frames in Hilbert space. 

I am grateful to M.~Rudelson for many important discussions.
The research would not be accomplished without help and 
encouragement of my wife Lilya.

\section{Preliminaries}                                                \label{SecPrelim}

We denote by $c, c_1, c_2$ absolute constants, 
and by $C(t)$, $C_1(t)$, $C_2(t)$ constants which depend on the
parameter $t$ only. The values of these constants may differ
from line to line.
The canonical vectors in $\R^n$ are denoted by $e_j$.

A sequence of vectors $(x_j)$ in a Banach space is called 
{\em $K$-Hilbertian} if 
$\| \sum a_j x_j \|  \le  K ( \sum |a_j|^2 )^{1/2}$
for any finite set of scalars $(a_j)$.
Similarly, $(x_j)$ is called 
{\em $K$-Besselian} if 
$K \| \sum a_j x_j \|  \ge  ( \sum |a_j|^2 )^{1/2}$
for any finite set of scalars $(a_j)$.
Suppose we are given two sequences $(x_j)$ and $(y_j)$ 
in Banach spaces $X$ and $Y$ respectively.
The sequences $(x_j)$ and $(y_j)$ are called {\em $K$-equivalent}
if there exist constants $K_1$ and $K_2$ with $K_1 K_2  \le  K$
such that for any finite sequence of scalars $(a_j)$
$$
K_1^{-1} \| \sum a_j y_j \|_Y
\le  \| \sum a_j x_j \|_X
\le  K_2 \| \sum a_j y_j \|_Y.
$$
In other words, the linear operator 
$T : \span(x_j)  \ra  \span(y_j)$ defined as
$T x_j  =  y_j$ for all $j$
is a $K$-isomorphism: $\|T\| \|T^{-1}\|  \le  K$.

Here and in the next section we work in a Hilbert space $H$
whose scalar product is denoted by $\< \cdot, \cdot \> $, 
and the norm by $\| \cdot \|$.
First we observe that the Hilbert-Schmidt norm of an operator on $H$
can be computed on the elements
of certain decompositions of identity.

\begin{lemma}                                                             \label{hs}
  Let $\id  =  \sum x_j \otimes x_j$
  be a decomposition on a Hilbert space $H$, 
  and $T : H \ra H$ be a linear operator. 
  Then 
  $$
  \|T\|_\HS^2  =  \sum \|T x_j\|^2.
  $$
\end{lemma}

\proof 
It is enough to write 
$$
\sum T x_j  \otimes  T x_j
=  T T^*
= \sum T e_j  \otimes  T e_j
$$
and to take traces.
\endproof

\noindent As an immediate consequence we have 

\begin{corollary}                                                             \label{squares}
  Let $\id  =  \sum x_j \otimes x_j$
  be a decomposition on a Hilbert space $H$. 
  Then 
  $$
  \sum \|x_j\|^2  =  \dim H.
  $$
\end{corollary}

\begin{lemma}
  Let $\id  =  \sum x_j \otimes x_j$
  be a decomposition on a Hilbert space.
  Then the system $(x_j)$ is $1$-Hilbertian. 
\end{lemma}

\proof
Notice that for every vector $x$
$$
\|x\|^2  =  \ll x, x \rr  
  =  \Big\ll \sum \ll x_j, x \rr x_j, x \Big\rr
  =  \sum | \ll x_j, x \rr |^2.
$$
Thus $\| \sum x_j \otimes e_j \| = 1$, 
and by duality $\| \sum e_j \otimes x_j \| = 1$.
This yields that $(x_j)$ is $1$-Hilbertian.
\endproof

The starting point of this paper is the 
principle of restricted invertibility proved 
by J.~Bourgain and L.~Tzafriri \cite{B-Tz}.

\begin{theorem} (J.~Bourgain, L.~Tzafriri).                          \label{BTz}  
  Let $T$ be a linear operator in $l_2^n$
  for which $\|T e_j\| = 1$, $j = 1, \ldots, n$. 
  Then there exists a subset $\s  \subset  \{1, \ldots, n\}$ 
  of cardinality $|\s|  \ge  c_1 n / \|T\|^2$ such that
  $$
  \Big\| \sum_{j \in \s} a_j T e_j \Big\|
    \ge  c \Big( \sum_{j \in \s} |a_j|^2 \Big)^{1/2}
  $$
  for any choice of scalars $(a_j)$.
\end{theorem}

The invertibility principle will be used together with
the following restriction theorem. It can easily be recovered 
from A.~Kashin's and L.~Tzafriri's paper \cite{K-Tz}, 
see the proofs of Theorem 1 and Corollary 2 there. 

\begin{theorem} (A.~Kashin, L.~Tzafriri).                                  \label{Lunin}
  Let $A$ be a norm-one linear operator in $l_2^m$.
  Fix a number $\l$ with $1/m  \le  \l  \le 1$.
  Then there exists a subset $\nu  \subset \{1, \ldots, m\}$
  of cardinality $|\nu|  \ge  \l m / 4$
  such that 
  $$
  \| P_\nu A \|  \le  c \left( \sqrt{\l} + \frac{\|A\|_\HS}{\sqrt{m}} \right).
  $$
Here $P_\nu$ denotes the coordinate projection onto $\R^\nu$.
\end{theorem}

Now we introduce an elementary procedure of splitting a sequence.
Given a sequence $(x_j)$ in $H$, let $(y_k)$ be any sequence of vectors in $H$ 
such that for every $j = 1, 2, \ldots$

\begin{itemize}
\item{the vectors $y_k$, $k \in \s_j$, are multiples of the vector $x_j$;}
\item{$\sum_{k \in \s_j} \|y_k\|^2  =  \|x_j\|^2$.}
\end{itemize}

\noindent Then we say that $(y_k)$ is the {\em splitted sequence} $(x_j)$.
Splitting allows us to make the norms of the vectors almost equal. Still the key 
property of a sequence, being $h$-Hilbertian, is preserved by splitting.

\section{The main result}                                               \label{SecMain}

In this section we prove an extraction theorem which is a core of the paper.

\begin{theorem}                                                          \label{main}
  Let $\id  =  \sum x_j \otimes x_j$
  be a decomposition on $l_2^n$,
  and $T$ be a norm-one linear operator. 
  Then for any $\e > 0$ there exists a set of indices $\s$
  of cardinality $|\s|  \ge  (1 - \e) \|T\|_\HS^2$ such that
  
  (i) the system $(T x_j)_{j \in \s}$ is $C(\e)$-equivalent
      to an orthogonal basis in $l_2^\s$;

  (ii) $\|T x_j\|  \ge  c \sqrt{\e} \frac{\|T\|_\HS}{\sqrt{n}} \|x_j\|$
      for all $j \in \s$.
\end{theorem} 

\proof
Put $h  =  \|T\|_\HS^2$.
By an approximation one can assume that the system $(x_j)$ is finite,
so we enumerate it as $(x_j)_{j \le m}$.
Denote $y_j  =  T x_j$ for all $j$. 
Splitting the system $(x_j)_{j \le m}$ we can assume that 
$$
0.9 \sqrt{ \frac{h}{m} }  \le  \|y_j\|  \le  1.1 \sqrt{ \frac{h}{m} }
\ \ \ \ \mbox{for all $j \le m$.}
$$
Let $\d  =  \e / 3$.
Consider the set of indices satisfying (ii), i.e.
$$
\t  =  \Big\{ j \le m :  
          \|y_j\|  \ge  0.9 \sqrt{\d} \sqrt{\frac{h}{n}} \|x_j\|  \Big\}.
$$
We claim that 
\begin{equation}                                            \label{sizeoftau}
  |\t|  \ge  (1 - \d) m.
\end{equation}
Indeed, by Corollary \ref{squares}
\begin{eqnarray*}
n  =  \sum_{j \le m} \|x_j\|^2  
    &\ge&  \sum_{j \in \t^c} \|x_j\|^2                                     \\
    &\ge&  \sum_{j \in \t^c} 
                  \left( \frac{1}{0.9 \sqrt{\d}} \sqrt{\frac{n}{h}} 
                  \|y_j\| \right)^2                                        \\
    &\ge&  |\t^c| \left( \frac{1}{0.9 \sqrt{\d}} \sqrt{\frac{n}{h}} 
                           0.9 \sqrt{\frac{h}{m}} \right)^2                
      =    |\t^c| \cdot \frac{n}{\d m}.
\end{eqnarray*}
Thus $|\t^c|  \le  \d m$, which proves (\ref{sizeoftau}). 

Now we have to find a further subset $\s  \subset \t$ 
of cardinality $|\s|  \ge  (1 - \e) h$,
such that the system $(y_j)_{j \in \s}$ is $C(\e)$-equivalent
to an orthogonal basis in $l_2^\s$.
The set $\s$ will be constructed by successive iterations. 
On the first step $\s = \emptyset$.
On each successive step, the remainder $h - |\s|$ 
will be reduced in a fixed proportion.
So, it is enough to prove the following.

\begin{lemma}
  Let $\s  \subset \t$ with $|\s|  <  (1 - \e)h$ be given, 
  and suppose the system $(y_j)_{j \in \s}$ is $K$-equivalent to an 
  orthogonal basis in Hilbert space.
  Then $\s$ can be extended in $\t$ to a subset $\s_1$ so that 
  
  (a) the system $(y_j)_{j \in \s_1}$ is $C(K, \e)$-equivalent to an
      orthogonal basis in Hilbert space;

  (b) for some absolute constant $\a < 1$ 
      $$
      h - |\s_1|  \le  \a (h - |\s|).                          \label{remainder}
      $$
\end{lemma}

\proof
Let $P$ be the orthonormal projection in $l_2^n$
onto $l_2^n  \ominus  \span(y_j)_{j \in \s}$
(at the first step $P = \id$).
First observe that by Lemma \ref{hs}
\begin{eqnarray*}
\sum_{j \le m} \|P y_j\|^2
  & = &  \sum_{j \le m} \|P T x_j\|^2
    =    \|PT\|_\HS^2                                          \\
  & = &  h  -  \| (\id - P) T \|_\HS^2                         \\
  &\ge&  h  -  \|\id - P\|_\HS^2 \|T\|^2                       \\
  & = &  h  -  |\s|.
\end{eqnarray*}
Using (\ref{sizeoftau}), we get
\begin{eqnarray*}
\sum_{j  \in  \t \setminus \s} \|P y_j\|^2  
  & = &  \sum_{j \in \t} \|P y_j\|^2                                       \\
  &\ge&  \sum_{j \le m}  \|P y_j\|^2  
            -  |\t^c| \cdot \Big( 1.1 \sqrt{h / m} \Big)^2                 \\
  &\ge&  h - |\s| - 1.21 \d h                                              \\
  &\ge&  (1 - 2 \d) h - |\s|   =:  h_0.
\end{eqnarray*}
Note that $h_0$ is comparable with $h$. 
Indeed, since $|\s|  <  (1 - \e) h  =  (1 - 3 \d) h$, we have
\begin{equation}
  h_0  \ge  \d h.                                              \label{hzeroh}
\end{equation}
Now we can split the system $(y_j)_{j  \in  \t \setminus \s}$
so that the resulting system $(y'_j)_{j \le M}$ satisfies
\begin{equation}
  \|P y'_j\|  \ge  0.9 \sqrt{\frac{h_0}{M}}                   
     \ \ \ \ \mbox{for all $j \le M$}.                                \label{Pydash}
\end{equation}

We are going to apply Kashin-Tzafriri's extraction result, Theorem \ref{Lunin}.
In the dual setting it can be reformulated as follows.

\begin{itemize}
\item{Let $(x_j)_{j \le m}$ be a $1$-Hilbertian system in a Hilbert space, 
  and put $\sum\|x_j\|^2  =  h$. 
  Fix a number $\l$ with $1/m  \le  \l  \le 1$.
  Then there exists a subset $\nu  \subset \{1, \ldots, m\}$
  of cardinality $|\nu|  \ge  \l m / 4$
  such that setting 
  $K  =  \Big( \sqrt{\l} + \sqrt{h / m} \Big)^{-1}$,
  $$
  \mbox{the system $(K x_j)_{j \in \s}$ is $c$-Hilbertian}.
  $$}
\end{itemize}

\noindent Apply this to the system $(y'_j)_{j \le M}$ 
which is $1$-Hilbertian and 
$$
\sum_{j \le M} \|y'_j\|^2  
  \le  \sum_{j \le m} \|y_j\|^2  
   =   h.
$$
With $\l  =  4 h / M$, we obtain a subset 
$\nu  \subset \{1, \ldots, M\}$
of cardinality $|\nu|  \ge  h$
such that 
\begin{equation}                                                    \label{ishilb}
\mbox{the system $\Big( \sqrt{\frac{M}{h}} y'_j \Big)_{j \in \nu}$ 
      is $c$-Hilbertian}
\end{equation}
(notice that we could make $M$ large enough to have $\l  \le  1$
as required in Kashin-Tzafriri's Theorem).
Therefore, the system 
$\Big( \sqrt{\frac{M}{h}} P y'_j \Big)_{j \in \nu}$ 
is $c$-Hilbertian, too.
Moreover, by (\ref{Pydash}) and (\ref{hzeroh})
$$
\Big\| \sqrt{\frac{M}{h}} P y'_j \Big\|   
  \ge  0.9 \sqrt{\frac{h_0}{h}}
  \ge  0.9 \sqrt{\d}  
  \ \ \ \ \mbox{for all $j \in \nu$}.
$$

At this point we use the original invertibility principle
of J.~Bourgain and L.~Tzafriri \cite{B-Tz}, 
which can be reformulated as follows.

\begin{itemize}
\item{Let $(x_j)_{j \le n}$ be an $H$-Hilbertian system in $l_2$ 
      and $\|x_j\|  \ge  \a$ for all $j$.
      Then there exists a subset $\rho  \subset \{1, \ldots, n\}$
      of cardinality $|\rho|  \ge  c (\a / H)^2 n$ such that
      the system $(x_j)_{j \in \rho}$ is $(c_1 / \a)$-Besselian.}
\end{itemize}

\noindent Apply this to the system 
$\Big( \sqrt{\frac{M}{h}} P y'_j \Big)_{j \in \nu}$.
There exists a subset $\rho'  \subset  \nu$ of cardinality 
$|\rho'|  \ge  c (0.9 \sqrt{\d})^2 h  =  c \d h$
such that
\begin{equation}                                               \label{isbess}
\mbox{the system $(\sqrt{\frac{M}{h}} P y'_j)_{j \in \rho'}$
      is $c_1 / \sqrt{\d}$-Besselian.}
\end{equation}
Recall that each vector $y'_j$ with $j \in \rho'$ 
is a multiple of some vector $y_{k(j)}$ 
with $k(j)  \in  \t \setminus \s$.
By (\ref{isbess}), these vectors $y_{k(j)}$ must be linearly independent. 
In particular, the correspondence $j \mapsto k(j)$ is one-to-one.
Consider the subset $\rho  \subset  \t \setminus \s$
consisting of vectors
$$
\rho  =  \{ k(j) : j \in \rho' \}.
$$
Now put $\s_1  =  \s \cup \rho$.

We see that (b) is satisfied with $\a  =  1 - c \d / 2$:
$$
h - |\s_1|  \le  h - |\s| - c \d h  \le  \a (h - |\s|).
$$

The reason why (a) holds is that 
the system $(y_j)_{j \in \rho}$ is well equivalent to an
orthogonal basis, 
and the spans of $(y_j)_{j \in \s}$ and $(y_j)_{j \in \rho}$
are well disjointed.
To implement this idea, note that by the preseding observations
(\ref{ishilb}) and (\ref{isbess}) yield that there exist constants 
$(\l_j)_{j \in \rho}$ such that

\begin{itemize}
\item{the system $(\l_j   y_j)_{j \in \rho}$ is $c$-Hilbertian,}
\item{the system $(\l_j P y_j)_{j \in \rho}$ 
                               is $c_1 / \sqrt{\d}$-Besselian.}
\end{itemize}

\noindent Now it is easy to complete the proof.
Fix any scalars $(a_j)_{j \in \rho}$.
Then defining $\l_j  =  y_j / \|y_j\|$ for $j \in \s$
we have
\begin{eqnarray*}
\Big( \sum_{j \in \s \cup \rho} |a_j|^2 \Big)^{1/2} 
  & = &  \Big( \sum_{j \in \s}  |a_j|^2 
             + \sum_{j \in \rho}|a_j|^2 \Big)^{1/2}                   \\
  &\le&  (K + c_1 / \sqrt{\d}) \left( 
         \Big\| \sum_{j \in \s}   a_j \l_j   y_j \Big\|^2
       + \Big\| \sum_{j \in \rho} a_j \l_j P y_j \Big\|^2
                               \right)^{1/2}                            \\
  & = &  (K + c_1 / \sqrt{\d}) 
         \Big\| \sum_{j \in \s}   a_j \l_j   y_j 
             +  \sum_{j \in \rho} a_j \l_j P y_j \Big\|       
            \\ \mbox{by orthogonality}                       \\
  &\le&  (K + c_1 / \sqrt{\d}) 
         \Big\| \sum_{j \in \s \cup \rho} a_j \l_j y_j \Big\|,
\end{eqnarray*}
because $P y_j = y_j$ for $j \in \s$.
This shows that the system $(\l_j y_j)_{j \in \s_1}$ 
is $(K + c_1 / \sqrt{\d})$-Besselian.
Next, 
\begin{eqnarray*}
\Big\| \sum_{j \in \s_1} a_j \l_j y_j \Big\|
  &\le&  \Big\| \sum_{j \in \s}   a_j \l_j y_j \Big\|
       + \Big\| \sum_{j \in \rho} a_j \l_j y_j \Big\|               \\
  &\le&  K \Big( \sum_{j \in \s}   |a_j|^2 \Big)^{1/2}        
       + c \Big( \sum_{j \in \rho} |a_j|^2 \Big)^{1/2}              \\
  &\le&  \sqrt{2} (K + c) 
           \Big( \sum_{j \in \s \cup \rho} |a_j|^2 \Big)^{1/2},
\end{eqnarray*}
so the system $(\l_j y_j)_{j \in \s_1}$ is $\sqrt{2}(K + c)$-Hilbertian as well. 
This establishes (a) of the Lemma
and completes the proof of Theorem \ref{main}.
\endproof

Let us rewrite Theorem \ref{main} in a different form,
which is useful in applications like Dvoretzky-Rogers type lemmas.

\begin{theorem}                                                            \label{mainc}
  Let $\id  =  \sum x_j \otimes x_j$
  be a decomposition on $l_2^n$,
  and $T$ be a norm-one linear operator, $\|T\|_\HS^2  =  h$.
  Then for any integer $\k  <  h$ there exists a set of indices $\s$
  with $|\s| = \k$ such that
  
  (i) the system $(T x_j)_{j \in \s}$ is $C(\k / h)$-equivalent
      to an orthogonal basis in $l_2^\s$;

  (ii) $\|T x_j\|  \ge  c \sqrt{\frac{h - \k}{n}} \|x_j\|$
      for all $j \in \s$.
\end{theorem}

\noindent It is sometimes more useful to get a lower bound 
for $\< x_j, T x_j \> $ rather than for $\|T x_j\|$.

\begin{proposition}                                                   \label{mainscalar}
  In Theorem \ref{main} statement (ii) can be replaced by 
  
  (ii') $| \< x_j, T x_j \> |  
          \ge  c \e \frac{|\trace T|}{n} \|x_j\|^2$
          for $j \in \s$.
\end{proposition}

\proof
Notice that $\sum  \< x_j, T x_j \>  =  \trace T$.
Therefore, splitting our system $(x_j)_{j \le m}$
we can assume that
\begin{equation}                                                       \label{tracem}
  | \< x_j, T x_j \> |  
    \ge  0.9 \frac{|\trace T|}{m}
    \ \ \ \ \mbox{for all $j$}.
\end{equation}

Let us examine the proof of Theorem \ref{main}.
The set $\tau$ was responsible for the lower bound of $\|T x_j\|$.
So, we replace $\tau$ by 
$$
\t'  =  \Big\{ j \le m :
       | \< x_j, T x_j \> |  
       \ge  (\e / 5) \frac{|\trace T|}{n} \|x_j\|^2
        \Big\},
$$
and all we have to check is
\begin{equation}                                                   \label{tdash}
  |\t'|  \ge  (1 - \e/3) m.
\end{equation}
By (\ref{tracem})
$$
\t'  \supset  \rho  
  =  \Big\{ j \le m : \|x_j\|^2  \le  (3 / \e) \frac{n}{m} \Big\}.
$$
Since $\sum_{j \le m} \|x_j\|^2  =  n$, 
we have $|\rho|  \ge  (1 - \e/3) m$.
This verifies (\ref{tdash}) and allows us 
to finish the proof as in Theorem \ref{main}.
\endproof

\section{Principle of restricted invertibility}                    \label{SecInvert}

Our first application stems from the look at Theorem \ref{main}
as an extension of the "principle of restricted invertibility",
Theorem \ref{BTz}, proved by J.~Bourgain and L.~Tzafriri.
Indeed, for the coordinate decomposition 
$id = \sum e_j \otimes e_j$ we get

\begin{corollary}                                                \label{BTzgeneralized}
  Let $T$ be a norm-one linear operator in $l_2$.
  Then for any $\e > 0$ there exists a subset $\s \subset \{1, 2, \ldots\}$
  of cardinality $|\s|  \ge  (1 - \e) \|T\|_\HS^2$ such that
  the sequence $(T e_j)_{j \in \s}$ is $C(\e)$-equivalent to an orthogonal
  basis in Hilbert space.
\end{corollary}

This theorem generalizes the invertibility principle in two ways. 
First, instead of requiring that {\em all} norms $\|T e_j\|$ 
be large, we can assume only largeness of their average, 
which is the Hilbert-Schmidt norm of $T$.
This makes the result independent of the dimension $n$ of the space,
$\|T\|_\HS^2$ being a natural substitute for the dimension $n$.

The second improvement is that we obtain the subset $\s$ with 
the largest possible cardinality.
Corollary \ref{BTzgeneralized} allows us to get 
$|\s|  \ge  (1 - \e) n / \|T\|^2$ in the original invertibility
principle for any $0 < \e < 1$
(while the Bourgain-Tzafriri's argument proves only the 
existence of such $\e$).
In some applications one really needs almost full percentage;
particularly, this is important in estimating the distance to the 
cube, see \cite{Sz-T}.

Notice that the infinite-dimensional analogs of Theorem \ref{main}
and Corollary~\ref{BTzgeneralized} hold, too.

\begin{proposition}
  Let $\id  =  \sum x_j \otimes x_j$ be a decomposition 
  on a Hilbert space, and $T$ be a linear operator which is 
  not Hilbert-Schmidt. 
  Then for any $\e > 0$ there exists an infinite subset $\s$
  such that the sequence $(T x_j)_{j \in \s}$ 
  is $(1 + \e)$-equivalent to an orthogonal basis in Hilbert space. 
\end{proposition}
  
\proof
The subset $\s$ is constructed by a standard induction argument, 
modulo the following claim:

\begin{itemize}
\item{For any finite dimensional subspace $E$
$$
\sup_j \dist (T x_j / \|T x_j\|, E)  =  1.
$$
}
\end{itemize}

Assume the contrary. 
Denoting the orthogonal projection in $l_2$ onto $E$ by $P$, we would have 
$$
\inf_j  \Big\| P (T x_j / \|T x_j\|) \Big\|  =  \d  >  0,
$$
that is 
$$
\|P T x_j\|  \ge  \d \|T x_j\| 
  \ \ \ \ \mbox{for all $j$.}
$$
By Lemma \ref{hs}, this yields
$$
\|PT\|_\HS  \ge  \d \|T\|_\HS  =  \infty.
$$
But the operator $PT$ has finite rank, 
thus $\|PT\|_\HS$ must be finite. 
This contradiction completes the proof. 
\endproof

An infinite dimensional analogue of Bourgain-Tzafriri's Theorem \ref{BTz}
says that, given a linear operator $T$ in $l_2$ with 
$\|T e_j\|  =  1$, $j = 1, 2, \ldots$, 
there exists a subset $\s$ of $\{1, 2, \ldots \}$
with upper density $\dens \s  \ge  c / \|T\|^2$, 
such that the sequence $(T e_j)_{j \in \s}$
is $c$-Besselian \cite{B-Tz}.

As we loose the normalizing condition $\|T e_j\| = 1$, nothing can be said
in general about the density of $\s$.
Indeed, let $(y_k)$ be the result of a splitting of a canonical basis
in $l_2$, so that the sets $\s_j$ from the definition of splitting
satisfy $|\s_j|  \ra  \infty$ as $j  \ra  \infty$.
There exist a norm-one operator $T$ in $l_2$ such that 
$\|T\|_\HS  =  \infty$ and
$T e_k  =  y_k$, $k = 1, 2, \ldots$.
However, each term $y_k$ in the sequence $(y_k)$ is repeated $|\s_j|$
times. Therefore, the upper density of any subset $\s$ satisfying
the conclusion of Corollary \ref{BTzgeneralized} must be zero.

Similarly, in some cases $\s$ must be a sparse set with respect 
to the dimension.
More precisely, in general the sequence $(T x_j)_{j \in \s}$ 
spans a subspace of infinite codimension.
This follows easily from a result of P.~Casazza and O.~Christensen
discussed in Section \ref{SecFrames}.

\section{Contact points and Dvoretzky-Rogers lemmas}                       \label{SecDR}

Here we apply Theorem \ref{main} to John decompositions. 
Let $X = (\R^n, \| \cdot\|)$ be a Banach space. 
The ellipsoid of maximal volume contained in $B(X)$ is unique 
and it is called {\em maximal volume ellipsoid} of $X$.
Suppose the euclidean structure on $X$ is chosen so that
the maximal volume ellipsoid coincides with the unit euclidean ball $D_n$.
Let us write a John decomposition on $X$:
\begin{equation}                                               \label{Johndecomposition}
  \id  =  \sum_{j = 1}^m  x_j \otimes x_j,
\end{equation}
where $x_j / \|x_j\|_X$ are some contact points of $B(X)$ with the John's 
ellipsoid (see \cite{T-J}, \S 15.3).
John decompositions can be considered as a subclass
in the class of all decompositions of type $\id  =  \sum x_j \otimes x_j$.
Conversely, each decomposition (\ref{Johndecomposition})
is a John decomposition for a suitable Banach space 
$X = (\R^n, \| \cdot\|)$, whose maximal volume ellipsoid is the unit
euclidean ball $D_n$.
Actually, the norm on $X$ is given by 
$\|x\|_X  =  \max_{j \le m} | \< x, \frac{x_j}{\|x_j\|_X} \> |$.
This result goes back to F.~John \cite{J}, 
although other proofs were found recently by 
K.~Ball \cite{B2} and A.~Giannopoulos and V.~Milman \cite{G-M}.
Therefore, working with contact points instead of decompositions
$\id  =  \sum x_j \otimes x_j$ we do not lose generality.

Recasting Theorem \ref{main} in this light, we have

\begin{corollary} 
  Let $X$ be an $n$-dimensional Banach space whose maximal volume ellipsoid
  is the unit euclidean ball. 
  Let $T$ be a linear operator with $\|T\|_{2 \ra 2}  \le  1$.
  Then for any $\e > 0$ there are contact points $x_1, \ldots, x_k$
  with $k  \ge  (1 - \e) \|T\|_\HS^2$
  such that the system $(Tx_j)_{j \le k}$ is $C(\e)$-equivalent 
  in $l_2$-norm to the canonical basis of $l_2^k$.
\end{corollary}

\noindent Moreover, the norms $\|T x_j\|_2$ are well bounded below:
$$
\|T x_j\|  \ge  C_1(\e) \frac{\|T\|_\HS}{\sqrt{n}}
\ \ \ \ \mbox{for all $j \le k$.}
$$

$T$ being an orthogonal projection, this result leads us to 
a new Dvoretzky-Rogers type lemma , which we will discuss now. 
Suppose $X$ is an $n$-dimensional Banach space whose John's 
ellipsoid is the unit euclidean ball. 
Let $Z$ be a $k$-dimensional subspace of $X$.
Dvoretzky-Rogers Lemma states that, given a positive integer
$\k < k$,
there is an orthonormal system $(z_j)_{j  \le  \k}$ in $Z$ such that
$$
\|z_j\|_X  \ge  \sqrt{\frac{k - \k + 1}{n}}
  \ \ \ \ \mbox{for all $j  \le  \k$}.
$$
Let us sketch the proof. By induction, it is enough to find {\em one}
vector $z$ in $Z$ such that $\|z\|_2 = 1$ and 
$\|z\|_X  \ge  \sqrt{\frac{k}{n}}$
(then substitute $Z$ by $Z  \ominus  \span(z)$ and repeat the argument).
By duality, this is equivalent to finding a functional 
$x^*  \in  B(X^*)$ with $\|P x^*\|_2  \ge  \sqrt{\frac{k}{n}}$,
where $P$ is the orthogonal projection onto $Z$.
Let $\id_{X^*}  =  \sum \l_j x_j^* \otimes x_j^*$
be a John decomposition in $X^*$, that is $\sum \l_j  =  n$
and $x_j^*$ are contact points of $B(X^*)$.
Then $P  =  \sum \l_j P x_j^* \otimes P x_j^*$.
Taking the trace, we get $k  =  \sum \l_j \|P x_j^*\|_2^2$.
So, there is a $j$ such that $\|P x_j^*\|_2^2  \ge  \sqrt{\frac{k}{n}}$.
This completes the proof.

However, this argument, as well as other known proofs of the 
Dvoretzky-Rogers Lemma, only establish the existence of the vectors $z_j$.
In contrast to that, the argument based on Theorem \ref{main}
provides information about their position.

\begin{theorem}                                                           \label{DR}
  Let $X$ be an $n$-dimensional Banach space 
  whose maximal volume ellipsoid is the unit euclidean ball.
  Let $P$ be an orthogonal projection, $\rank P = k$.
  Then for any positive integer $\k  <  k$
  there are contact points $x_1, \ldots, x_\k$
  such that setting $z_j  =  P x_j / \|P x_j\|_2$
  we have
  
  (i)  the system $(z_j)$ is $C(\k / k)$-equivalent in $l_2$-norm 
       to the canonical basis of $l_2^\k$;

  (ii) $\|z_j\|_X  \ge  c \sqrt{\frac{k - \k}{n}}$
       for all $j$.
\end{theorem}

\proof
Note that $\|P\|_\HS  =  \sqrt{k}$ and apply Theorem \ref{mainc}
to a John decomposition on $X$. 
This gives us contact points $x_1, \ldots, x_\k$ such that
(i) is satisfied, and
$$
\|P x_j\|_2  \ge  c \sqrt{\frac{k - \k}{n}}
  \ \ \ \ \mbox{for every $j \le m$}.
$$
It remains to note that for every $j \le m$
$$
\|P x_j\|_X  \ge  \frac{\< P x_j, x_j \> }{\|x_j\|_{X^*}}
  =  \< P x_j, x_j \>  =  \|P x_j\|_2^2,
$$
hence
$$
\|z_j\|_X  \ge  \|P x_j\|_2
  \ \ \ \ \mbox{for every $j \le m$}.
$$
This completes the proof.
\endproof

A natural question is whether Theorem \ref{DR} 
can be extended for arbitrary operator $T$, 
$\|T\|_{2 \ra 2} = 1$, 
with $k$ substituted by $\|T\|_\HS^2$. 
The answer is negative. 
Indeed, let $n  =  2^m$ for a positive integer $m$,
and denote by $W_m$ the Walsh $n \times n$ matrix.
Consider the operator $T  =  n^{-1/2} W_m$
acting in the space $X  =  l_\infty^n$.
All contact points $(x_j)$ of $X$ are simply the
coordinate vectors. However, denoting $z_j = T x_j / \|T x_j\|_2$
we have for any $j$
$$
\|z_j\|_X  =  \|T x_j\|_X / \|T x_j\|_2  =  n^{-1/2}.
$$
This shows that (ii) in Theorem \ref{DR} would fail 
for the operator $T$.

Still, a lower bound for $\| \cdot \|_X$-norm exists
and is equivalent to $\frac{1}{n} |\trace{T}|$.

\begin{proposition}
  Let $X$ be an $n$-dimensional Banach space
  whose maximal volume ellipsoid is the unit euclidean ball.
  Let $T$ be an operator with $\|T\|_{2 \ra 2} \le 1$.
  Then for any $\e > 0$
  there are $k  >  (1 - \e) n$ contact points $x_1, \ldots, x_\k$
  such that 
  
  (i)  the system $(T x_j)$ is $C(\e)$-equivalent in $l_2$-norm 
       to the canonical basis of $l_2^\k$;

  (ii) $\|T x_j\|_X  \ge  c \e \frac{|\trace T|}{n}$.
\end{proposition}

\proof
The argument is similar to that of Theorem \ref{DR};  
one only uses Proposition \ref{mainscalar}
instead of Theorem \ref{main}.
\endproof

The estimate in (ii) is essentially sharp. 
Indeed, for any positive integers $k \le n$
one can construct an orthogonal projection $P$
in $\R^n$ with $\rank P = k$, and such that
$$
\|P e_j\|_2  =  \sqrt{\frac{k}{n}}
\ \ \ \ \mbox{for all $j \le n$}.
$$
Notice that 
$\|P\|_{1 \ra 2}       =  \sqrt{\frac{k}{n}}$, 
therefore 
$\|P\|_{2 \ra \infty}  =  \sqrt{\frac{k}{n}}$.
Thus 
$\|P\|_{1 \ra \infty}  \le  \frac{k}{n}$.
Now consider the space $X  =  l_\infty^n$
and the operator $T = P$ on it.
The only contact points $(x_j)$ of $X$ are the coordinate 
vectors $(e_j)$. Then for all $j$
$$
\|T x_j\|_X  =  \|P e_j\|_\infty  
  \le  \frac{k}{n}  
   =  \frac{|\trace T|}{n}.
$$
This shows that the lower bound in (ii) is essentially sharp.

Finally, there is a class of operators for which 
Theorem \ref{DR} itself can be extended: selfadjoint operators.
So, the desired general result for arbitrary operator can be obtained
if we allow a suitable rotation of $(z_j)$.
There is always a unitary operator (coming from the polar
decomposition of $T$) which sends the vectors $(z_j)$
to vectors of a good $\| \cdot \|_X$-norm.

\begin{theorem}                                                    
  Let $X$ be an $n$-dimensional Banach space 
  whose maximal volume ellipsoid is the unit euclidean ball.
  Let $T$ be an operator with $\|T\|_{2 \ra 2} \le 1$, 
  and put $\|T\|_\HS^2  =  h$.
  Then for any positive integer $\k  <  h$
  there are contact points $x_1, \ldots, x_\k$
  such that setting $z_j  =  |T| x_j / \|T x_j\|_2$
  we have
  
  (i)  the system $(z_j)$ is $C(\k / h)$-equivalent in $l_2$-norm 
       to the canonical basis of $l_2^\k$;

  (ii) $\|z_j\|_X  \ge  c \sqrt{\frac{h - \k}{n}}$
       for all $j$.
\end{theorem}

\proof
Let $T = U |T|$ be the polar decomposition of $T$, 
where $|T| = (T^* T)^{1/2}$ is a positive selfadjoint operator
and $U$ is a partial isometry on $l_2^n$.
Apply Theorem \ref{mainc} to the operator $|T|$
and a John decomposition on $X$. 
As before, this gives us contact points $x_1, \ldots, x_\k$ 
such that (i) is satisfied, and
$$
\Big\| |T| x_j \Big\|_2  \ge  c \sqrt{\frac{h - \k}{n}}
  \ \ \ \ \mbox{for every $j \le m$}.
$$
From diagonalization of $|T|$ it follows that
$\Big\| |T|^{1/2} \Big\|_{2 \ra 2}  =  \|T\|_{2 \ra 2}^{1/2}  \le  1$.  
Therefore we can bound for every $j \le m$
\begin{eqnarray*}
\Big\| |T| x_j \Big\|_X  
 &\ge&  \frac{\< |T| x_j, x_j \> }{\|x_j\|_{X^*}}
   =    \Big\| |T|^{1/2} x_j \Big\|_2^2                \\
 &\ge&  \Big\| |T| x_j \Big\|_2^2
   =    \Big\| T x_j \Big\|_2^2.
\end{eqnarray*}
Hence
$$
\|z_j\|_X  \ge  \|T x_j\|_2
  \ \ \ \ \mbox{for every $j \le m$}.
$$
This completes the proof.
\endproof

Now we turn to a particular case when $T$ is the identity operator, 
which also happens to be interesting. We clearly have 

\begin{corollary}                                            \label{contact}
  Let $X$ be an $n$-dimensional Banach space 
  whose maximal volume ellipsoid is the unit euclidean ball.
  Then for any $\e > 0$ there is a set of $k  >  (1 - \e) n$
  contact points which is $C(\e)$-equivalent to 
  the canonical vector basis of $l_2^k$.
\end{corollary}

\noindent This result is related to another variant of the 
classical Dvoretzky-Rogers Lemma, which establishes the existence
of contact points whose distance to a certain orthonormal basis
is controlled (\cite{T-J}, Theorem 15.7).
More precisely, there exist contact points $x_1, \ldots, x_n$
and an orthonormal basis $h_1, \ldots, h_n$ such that
$$
\|x_j - h_j\|_2  \le  2 \Big(1 - \sqrt{\frac{n - i + 1}{n}} \Big)
\ \ \ \ \mbox{for $i  \le  n$}.
$$
However, this estimate is to crude to assure that a fixed 
proportion of the system $(x_j)_{j \le n}$ is equivalent 
in $l_2$-norm to an orthonormal system. Such an isomorphism 
can be established only for $c \sqrt{n}$ contact points. 

Using this argument, it is proved that for $k = [\sqrt{n} / 4]$
there exist orthonormal vectors 
$x_1, \ldots, x_k$ in $(E, \|\cdot\|_2)$ on which all three norms 
$\|\cdot\|$, $\|\cdot\|_2$, and $\|\cdot\|_*$ 
differ by the factor $2$ at most (\cite{T-J}, p.127). 
It has been an open problem to make $k$ proportional to $n$.
By Corollary \ref{contact} we actually have $k  \ge  (1 - \e)n$
and make all three norms $\|\cdot\|$, $\|\cdot\|_2$, and $\|\cdot\|_*$
equal to $1$ on our sequence (paying however in exact orthogonality).

By duality, Corollary \ref{contact} holds also for the ellipsoid of minimal 
volume containing $B(X)$ instead of the maximal volume ellipsoid. 
This variant of Corollary \ref{contact} yields also a proportional 
Dvoretzky-Rogers factorization result from \cite{B-Sz}. 
Namely, given an $n$-dimensional Banach space $X$ and $\e > 0$, 
there is a $k > (1 - \e) n$ such that the formal identity 
$\id : l_2^k \ra l_\infty^k$
can be written as $\id = \a \cdot \b$ for some 
$\b : l_2^k \ra X$, $\a : X \ra l_\infty^k$, 
with $\|\a\| \|\b\|  \le  C(\e)$.
This can be obtained by duality from the factorization of the identity
on the contact points, $\id : l_1^k  \ra  X  \ra  l_2^k$, 
guaranteed by Corollary \ref{contact}.

A few comments about the dependence $C(\e)$ in Theorem \ref{main}.
It is a challenge to find the correct asympthotics. 
Indeed, by an argument of S.~Szarek and M.~Talagrand \cite{Sz-T}
the proportional Dvoretzky-Rogers factorization above
yields a non-trivial estimate on the distance from $X$ to $l_\infty^n$ --
a well known and hard problem in the local theory of normed spaces.
The factorization constant $C(\e)$ lies in the heart of the computation
of this distance.

The proof of Theorem \ref{main} guarantees that
$C(\e)  \le  \e^{c \log \e}$.
However, for the Dvoretzky-Rogers factorization constant 
$C_\DR (\e)$ much better bounds are found \cite{G}: 
$C_\DR (\e)  \le  c \e^{-1}$. 
Since $C_\DR (\e)  \le  C(\e)$ and 
$C_\DR (\e)  \ra  \infty$ as $\e \ra 0$  (see \cite{Sz-T}), 
we necessarily have $C(\e)  \ra  \infty$ as $\e \ra 0$.

However, $C(\e)  \ra  1$ as $\e \ra 1$.
This follows directly from 

\begin{lemma}
  Let a normalized sequence $(x_j)_{j \le n}$ in Hilbert space 
  be $M$-Hilbertian, and $\e > 0$.
  Then there is a subset $\s  \subset  \{1, \ldots, n\}$
  of cardinality $|\s|  \ge  C(M, \e) n$ such that
  the system $(x_j)_{j \in \s}$ is $(1 + \e)$-equivalent
  to the canonical basis of $l_2^\s$.
\end{lemma}

\proof
We can assume that the given Hilbert space is $l_2^n$.
Define a linear operator $T : l_2^n  \ra  l_2^n$ by 
$$
T e_j  =  x_j  
\ \ \ \ \mbox{for $j \le n$.}
$$
Let $A = T^*T - \id$.
Then the matrix of $A$ has zeros on the diagonal and 
$\|A\|  \le  M^2 + 1$.
Now, by a theorem of J.~Bourgain and L.~Tzafriri
(\cite{B-Tz} Theorem 1.6, see also \cite{K-Tz})
there is a subset $\s  \subset  \{1, \ldots, n\}$ 
of cardinality $|\s|  \ge  C(M, \e)$ such that
$$
\| P_\s A P_\s \|  <  \e.
$$
This shows that for any sequence of scalars $(a_j)_{j \le n}$
$$
\Big| \Big\ll
(T^*T - \id) \sum_{j \in \s} a_j e_j , \sum_{j \in \s} a_j e_j
\Big\rr \Big|
<  \e,
$$
hence
$$
\Big| \Big\ll
\sum_{j \in \s} a_j x_j , \sum_{j \in \s} a_j x_j
\Big\rr 
  -  \sum_{j \in s} |a_j|^2
\Big|
<  \e.
$$
This clearly finishes the proof.
\endproof

\section{Embeddings of the cube}                                     \label{SecCube}

In this section we apply Theorem \ref{main} to the study 
of embeddings of $l_\infty^k$ into finite dimensional spaces. 
N.~Alon and V.~Milman proved that if a given normalized sequence $(x_j)$ 
in a Banach space $X$ has small Rademacher average $\E \| \sum \e_j x_j\|_X$, 
then it must contain a large subsequence well equivalent 
to the canonical basis of $l_\infty^k$.
Later on, M.~Talagrand \cite{T} improved this result 
and simplified the argument.

\begin{theorem} (M.~Talagrand).  \                                   \label{Talagrand}
  Suppose we are given vectors $(x_j)_{j \le n}$ in a Banach space $X$
  with $\|x_j\|_X  \ge  1$.
  Set $M = \E \| \sum \e_j x_j \|_X$ 
  and $\w = \sup \Big\{ \sum |x^*(x_j)| : x^* \in B(X^*) \Big\}$.
  Then there exists a subset $\s  \subset  \{1, \ldots, n\}$ 
  of cardinality $|\s|  \ge  c n / \w$ such that
  $$
  \frac{1}{2} \max_{j \in \s} |a_j|
    \le  \Big\| \sum_{j \in \s} a_j x_j \Big\|_X
    \le   4 M \max_{j \in \s} |a_j|
  $$  
  for any choice of scalars $(a_j)$.
\end{theorem}

A few years earlier, M.~Rudelson obtained in \cite{R1} 
a "Gaussian" version of this theorem.
Recall that the $\ell$-norm of an operator $u : l_2^n  \ra  X$
is defined as 
$$
\ell(u)^2  =  \int \|ux\|_X^2 \, d \gamma_n(x),
$$
where $\gamma_n$ is the canonical Gaussian measure on $\R^n$.
We will sometimes write $\ell(X)$ instead of $\ell(\id_X)$.

Suppose $X$ is an $n$-dimensional Banach space whose maximal volume ellipsoid
is the unit euclidean ball.
Let $P$ be an orthogonal projection in $X$ 
and set $k  =  \rank P$, $a  =  k / n$.
The result of M.~Rudelson states that there is a subspace 
$Z  \subset  P(X)$  of dimension $m  \ge  C_1(a) \frac{\sqrt{n}}{\ell(P)}$
which is $C_2(a) \ell(P)$-isomorphic to $l_\infty^m$.

We will remove $\ell(P)$ from the estimate on the dimension.
Further, it will be shown that $Z$ is canonically spanned by the projections
of some contact points of $X$.
In particular, the norm on $Z$ is well equivalent to 
$\max_{j \le m} | \< z , x_j \> |$,
where $x_j$ are contact points.
This yields automatically that $Z$ is well complemented 
by the orthogonal projection. 
Moreover, the dependence on $a$ will be improved.

\begin{theorem}                                                    \label{GT}
  Let $X$ be an $n$-dimensional Banach space 
  whose maximal volume ellipsoid is $B(l_2^n)$.
  Let $P$ be an orthogonal projection, $\rank P = k$.
  Then there are contact points $(x_j)_{j  \le  m}$
  with $m  \ge  c_1 k / \sqrt{n}$ such that 
  $$
  \max_{j \le m} | \< x , x_j \> |    \le   \|x\|_X    
  \le   c \sqrt{\frac{n}{k}} \ell(P) \max_{j \le m} | \< x , x_j \> |
  $$
  for every $x$ in $Z  =  \span(P x_j)_{j \le m}$.
\end{theorem}

A particular case $k = n$ is also interesting: 
we get a sequence of $m  \ge  c_1 \sqrt{n}$ contact points 
which is $c \ell(X)$-equivalent to the canonical vector basis 
of $l_\infty^m$.
Let us prove this latter fact separately. 
By Corollary \ref{contact}, there exists a set of contact points 
$(x_j)_{j \le m'}$, $m'  \ge  n / 2$, 
which is $c$-equivalent in $l_2$-norm 
to the canonical basis of $l_2^{m'}$.
Let $(g_j)$ be independent standard Gaussian random variables.
To apply Talagrand's Theorem \ref{Talagrand}, we bound
$$
M   =   \E \| \sum \e_j x_j \|_X
   \le  c \E \| \sum g_j x_j \|_X
   \le  c c_2 \E \| \sum g_j e_j \|_X
$$
by the ideal property of the $\ell$-norm, see Lemma \ref{Slepian} below. 
Further, for any finite system of scalars $(a_j)$
\begin{eqnarray*}
\| \sum a_j x_j \|_X   
  &\le&   \| \sum a_j x_j \|_2                     \\
  &\le& c_2 \Big( \sum |a_j|^2 \Big)^{1/2}          
   \le  c_2 \sqrt{n} \max_j |a_j|.
\end{eqnarray*}
Thus
$$
\w  =  \sup \Big\{ \sum |x^*(x_j)| : x^* \in B(X^*) \Big\}
   \le c_2 \sqrt{n}.
$$
Now Talagrand's Theorem \ref{Talagrand} finishes the proof. 
\endproof

The proof of Theorem \ref{GT} is longer, but the main idea remains 
to combine results of Section \ref{SecDR} with Talagrand's theorem. 
First, we need to know what vectors canonically span a large
subspace $Z$ well isomorphic to $l_\infty^m$.
They happen to be multiples of some contact points $(x_j)_{j \le m}$.
More precisely, we have

\begin{proposition}                                                           \label{spancube}
  Under the assumptions of Theorem \ref{GT}, the system 
  $\Big( \frac{Px_j}{\|Px_j\|_2} \Big)_{j \le m}$
  is $c \sqrt{\frac{n}{k}} \ell(P)$-equivalent 
  to the canonical basis of $l_\infty^m$.
\end{proposition}

To prove this, let $(x_j)_{j \le k'}$, $k'  \ge  k / 2$,
be the contact points provided by Theorem \ref{DR}.
Put 
$$
z_j  =  \a \sqrt{\frac{n}{k}} \frac{Px_j}{\|Px_j\|_2}
$$
for an appropriate $\a > 0$ and all $j$.
Then 

(i)  the system $(z_j)_{j \le k'}$ is $(c \sqrt{\frac{n}{k}})$-equivalent
     in $l_2$-norm to the canonical basis of $l_2^{k'}$;

(ii) $\|z_j\|  \ge  1$ for all $j$.

\noindent For a future reference note that the proof of 
Theorem \ref{DR} gives also 
\begin{equation}                                                                 \label{Pxj}
  \|Px_j\|_2  \ge  c_1 \sqrt{\frac{n}{k}}.
\end{equation}

To apply Talagrand's Theorem \ref{Talagrand} to the system $(z_j)_{j \le k'}$,
recall the ideal property of the $\ell$-norm (cf. \cite{T-J}, \S 12)

\begin{lemma}                                                                   \label{Slepian}
  For any two linear operators 
  $A : l_2^n  \ra  X$ and $B : l_2^n  \ra  l_2^n$
  $$
  \ell(AB)  \le  \|B\| \ell(A).
  $$
\end{lemma}
Let $(h_j)_{j \le k'}$ be an orthonormal basis in $\span(z_j)$.
Now we bound 
\begin{eqnarray*}
M  & = &    \E \| \sum \e_j z_j \|_X  
    \le   c \E \| \sum  g_j z_j \|_X                          \\
   &\le&  c \sqrt{\frac{n}{k}} \E \| \sum g_j h_j \|_X         
          \ \ \ \ \mbox{by (i) and Lemma \ref{Slepian}}       \\
   & = &  c \sqrt{\frac{n}{k}} \ell(P(X))
     =    c \sqrt{\frac{n}{k}} \ell(P).
\end{eqnarray*}
Noting (ii) above, we apply Talagrand's Theorem \ref{Talagrand}.
This finishes the proof. 
\endproof

To obtain Theorem \ref{GT}, Talagrand's Theorem will be used more delicately.
Its proof in \cite{T} gives the following additional property.

\begin{lemma}                                                                 \label{Talagrandplus}
  In the situation of Theorem \ref{GT}, suppose $(x_j^*)_{j \le n}$
  are functionals in $X^*$ such that 
  $$
  x_j^* (x_j)  \ge  1  \ \ \ \ \mbox{and} \ \ \ \ 
  \|x_j^*\|_{X^*}  =  1      \ \ \ \ \mbox{for all $j \le n$.}
  $$
  Then the subset $\s$ can be found so that 
  $$
  \sum_{j  \in  \s \setminus \{i\}} |x_i^* (x_j)|  \le  1/2
  \ \ \ \ \mbox{for all $i \in \s$}.
  $$
\end{lemma}

Let us turn again to the proof of Proposition \ref{spancube}.
We applied Talagrand's Theorem to the system $(x_j)_{j \le k'}$.
Note that by (\ref{Pxj})
\begin{eqnarray*}
\< x_j, z_j \>  & = &  \a \sqrt{\frac{n}{k}} \frac{\< x_j, Px_j \> }{\|Px_j\|_2}  \\
  & = &   \a \sqrt{\frac{n}{k}} \|Px_j\|_2
   \ge  1   \ \ \ \ \mbox{for all $j$.}
\end{eqnarray*}
Therefore we obtain a subset $\s  \subset \{1, \ldots, k'\}$
of cardinality $|\s|  \ge  c_1 k / \sqrt{n}$ such that
\begin{itemize}
  \item{the system $(z_j)_{j \in \s}$ 
        is $c \sqrt{\frac{n}{k}} \ell(P)$-equivalent to the 
        canonical basis of $l_\infty^\s$;}
  \item{$\sum_{j  \in  \s \setminus \{i\}} |\< x_i , z_j\> |  \le  1/2$
        for all $j \in \s$.} 
\end{itemize}
Now fix an $x  =  \sum_{j \in \s} a_j z_j$,
and let $i \in \s$ be such that $|a_i|  =  \max_{j \in \s} |a_j|$.
Then
\begin{eqnarray*}
\max_{j \in \s} |\< x, x_j \> |  
  &\ge&  |\< x_i, x \> |
    =    \Big| \sum_{j \in \s} a_j \< x_i, z_j \> \Big|          \\
  &\ge&  |a_j| \Big( |\< x_i, z_i \> |  
           -  \sum_{j  \in  \s \setminus \{i\}} | \< x_i, z_j \> | \Big) \\
  &\ge&  \frac{1}{2} |a_i|
    =    \frac{1}{2} \max_{j \in \s} |a_j|                    \\
  &\ge&  \Big( c \sqrt{\frac{n}{k}} \ell(P) \Big)^{-1} \|x\|_X.
\end{eqnarray*}
This proves the second inequality in Theorem \ref{GT}, 
while the first one is obvious.
\endproof

Theorem \ref{GT} yields also that $Z$ is well complemented by 
the orthogonal projection. 

\begin{proposition}
  In the situation of Theorem \ref{GT}, 
  let $P_Z$ be the orthogonal projection in $X$ onto $Z$. Then 
  $$
  \|P_Z\|  \le  c \sqrt{\frac{n}{k}} \ell(P).
  $$
\end{proposition}

\proof 
For any $x \in X$
\begin{eqnarray*}
\|P_Z x\|_X  
  &\le&  c \sqrt{\frac{n}{k}} \ell(P) \max_{j \le m} | \< P_Z x, x_j \> | \\
  & = &  c \sqrt{\frac{n}{k}} \ell(P) \max_{j \le m} | \< x, x_j \> |     \\
  &\le&  c \sqrt{\frac{n}{k}} \ell(P) \|x\|_X.
\end{eqnarray*}
This completes the proof.
\endproof

Another consequence of Theorem \ref{GT} is a refined isomorphic 
characterization of spaces with large volume ratio.
Recall that the volume ratio of an $n$-dimensional Banach space $X$
is defined as
$$
\vr(X)  =  \min \left( \frac{\vol(B_X)}{\vol(E)} \right)^{1/n}
$$
over all ellipsoids $E$ contained in $B_X$, see \cite{L-M}, \cite{P}.
The maximal value of $\vr(X)$ among all $n$-dimensional spaces 
is of order $\sqrt{n}$, 
and the only space with maximal volume ratio is $l_\infty^n$
(K.~Ball \cite{B1}).
Later on, M.~Rudelson proved in \cite{R1} that
if $\vr(X)$ is proportional to the maximal volume ratio 
then $X$ has a subspace isomorphic to $l_\infty^m$ with 
the isomorphism constant of order $\log n$, 
and such that 
$m  \sim  \frac{\sqrt{n}}{\log n}$.
Using Theorem \ref{GT} in the M.~Rudelson's proof of this
result removes the parasitic factor $\log n$ from the estimate 
on the dimension. 

\begin{theorem}
  Let $a > 0$, and $X$ be an $n$-dimensional Banach space. 
  If $\vr(X)  \ge  a \sqrt{n}$ 
  then there exists a subspace $Z$ of dimension 
  $m  \ge  C_1(a) \sqrt{n}$ which is $C_2(a) \log n$-isomorphic 
  to $l_\infty^m$.
\end{theorem}

\section{Frames}                                                     \label{SecFrames}

The notion of frame goes back to the work of R.~Duffin and A.~Schaeffer
on nonharmonic Fourier series \cite{D-S}.
A sequence $(x_j)$ in a Hilbert space $H$ is called a {\em frame} 
if there exist positive numbers $A$ and $B$ such that
$$
A \|x\|^2  \le  \sum_j | \ll x, x_j \rr |^2  \le  B \|x\|^2
\ \ \ \ \mbox{for $x \in H$.}
$$
The number $(B/A)^{1/2}$ is called a {\em constant} of the frame.
We call $(x_j)$ a {\em tight frame} if $A = B = 1$.
For introduction to the theory of frames, its relation to wavelets and signal 
processing, see \cite{B-W}.
Geometric structure of frames is studied extensively in recent years, see
\cite{Ho}, \cite{A}, \cite{C-C1}, \cite{C-C2}, \cite{C2}, \cite{V}.

It is known by now that finite dimensional frames are essentially the 
same object as John decompositions. 
In the equivalent theory, it is sufficient to work only with tight frames,
because every frame with constant $M$ is $M$-equivalent to a tight frame
(cf. e.g. \cite{Ho}).
Further, one has the following equivalence between frames and 
John decompositions:
$$
\mbox{ $(x_j)$ is a tight frame in $H$ }   
  \iff   \id_H  =  \sum_j x_j \otimes x_j.
$$
This observation allows to interpret the results of 
Sections \S \ref{SecMain} and \S \ref{SecInvert} as statements about frames.
Theorem \ref{main} yields:

\begin{corollary}                                                    \label{frame}
  Let $(x_j)$ be a tight frame in Hilbert space $H$,
  and $T$ be a norm-one linear operator in $H$. 
  Then for any $\e > 0$ there exists a subset of indices $\s$
  of cardinality $|\s|  \ge  (1 - \e) \|T\|_\HS^2$ such that
  the system $(T x_j)_{j \in \s}$ is $C(\e)$-equivalent
  to an orthogonal basis in Hilbert space.
\end{corollary}
For operators $T$ which are not Hilbert-Schmidt this means 
that the subset $\s$ is infinite.

Clearly, Corollary \ref{frame} itself generalizes the invertibility
principle of J.~Bourgain and L.~Tzafriri.
When applied to the identity operator, it yields that every tight 
frame in $l_2^n$ has a subset of length $(1 - \e) n$ 
which is $C(\e)$-equivalent to an orthogonal basis in Hilbert space.
This result was proved in \cite{V} as a generalization of 
P.~Casazza's theorem \cite{C2}.

Notice that one necessarily has $C(\e)  \ra  \infty$ as $\e  \ra  1$,
as explained in Section \S \ref{SecDR}.
An infinite dimensional analog of this phenomenon holds, too.
A frame may not in general contain a complete subsequence equivalent
to an orthogonal basis. The counterexample was found 
by P.~Casazza and O.~Christensen in \cite{C-C2}, see also \cite{V}.


{\small
\begin{thebibliography}{S 99} 

\bibitem [A] {A}    A. Aldroubi, 
  {\em Portraits of frames},
  Proc. of the AMS 123 (1995), 1661--1668

\bibitem [B1] {B1}    K. Ball,
  {\em Volumes of sections of cubes and related problems},
  GAFA Seminar 87-88, Springer Lecture Notes in Math. 1376 (1989), 251--263

\bibitem [B2] {B2}    K. Ball,
  {\em Ellipsoids of maximal volume in convex bodies},
  Geom. Dedicata 41 (1992), 241--250

\bibitem [B-Sz] {B-Sz}    J. Bourgain, S. Szarek,
  {\em The Banach-Mazur distance to cube and the Dvoretzky-Rogers factorization},
  Israel J. Math. 62 (1988), 169--180

\bibitem [B-Tz] {B-Tz}    J. Bourgain,  L. Tzafriri,
  {\em Invertibility of "large" submatrices with applications 
       to the geometry of Banach spaces and harmonic analysis},
  Israel J. Math. 57 (1987), 137--224

\bibitem [B-W] {B-W}    J. Benedetto, D. Walnut, 
  {\em Gabor frames for $L^2$ and related spaces},
  in: Wawelets: mathematics and applications, ed. by J.~Benedetto, M.~Frazier, 
  CRC Press, 1994

\bibitem [C2] {C2}    P. G. Casazza, 
  {\em Local theory of frames and Schauder bases for Hilbert space},
  Illinois J. Math., to appear (1999) 

\bibitem [C-C1] {C-C1}    P. G. Casazza, O. Christensen,
  {\em Hilbert space frames containing a Riesz basis 
       and Banach spaces which have no subspace isomorphic to $c\sb 0$},
  J. Math. Anal. Appl. 202 (1996), 940--950  

\bibitem [C-C2] {C-C2}    P. G. Casazza, O. Christensen,
  {\em Frames containing a Riesz basis 
       and preservation of this property under perturbations},
  SIAM J. Math. Anal. 29 (1998), 266--278     

\bibitem [D-S] {D-S}    R. J. Duffin, A. C. Schaeffer,
  {\em A class of nonharmonic Fourier series},
  Trans. of the AMS 72 (1952), 341--366

\bibitem [G] {G}    A. Giannopoulos, 
  {\em A proportional Dvoretzky-Rogers factorization result},
  Proc. AMS 124 (1996), 233--241

\bibitem [G-M] {G-M}    A. Giannopoulos, V. Milman,
  {\em Extremal problems and isotropic positions of convex bodies},
  Preprint

\bibitem [G-P-T] {G-P-T}    A. Giannopoulos, I. Perissinaki, A. Tsolomitis,
  {\em John's theorem for an arbitrary pair of convex bodies},
  Preprint

\bibitem [Ho] {Ho}    J. R. Holub,
  {\em Pre-frame operators, Besselian frames, and near-Riesz bases
       in Hilbert spaces},
  Proc. of the AMS 122 (1994), 779--785

\bibitem [J] {J}    F. John, 
  {\em Extremum problems with inequalities as subsidiary conditions},
  Courant Anniversary Volume, Interscience, New York, 1948, 187--204

\bibitem [K-Tz] {K-Tz}    B. Kashin, L. Tzafriri, 
  {\em Some remarks on the restrictions of operators to coordinate subspaces},
  Preprint

\bibitem [L-M] {L-M}   J. Lindenstrauss, V. Milman,  
  {\em The local theory of normed spaces and its applications to convexity},
  Handbook of convex geometry, Vol. A, B, 1149--1220, 
  North-Holland, Amsterdam, 1993

\bibitem [P] {P}    G. Pisier,
  {\em The volume of convex bodies and Banach space geometry},
  Cambridge Tracts in Mathematics, 94, 
  Cambridge University Press, 1989

\bibitem [R1] {R1}    M. Rudelson,
  {\em Estimates of the weak distance between finite-dimensional Banach spaces},
  Israel J. Math. 89 (1995), 189--204

\bibitem [R2] {R2}    M. Rudelson,
  {\em Contact points of convex bodies},
  Israel J. of Math. 101 (1997), 93--124

\bibitem [Sz-T] {Sz-T}    S. Szarek, M. Talagrand,
  {\em An "isomorphic" version of the Sauer-Shelah lemma and the Banach-Mazur
       distance to the cube},
  GAFA Seminar 87-88, Springer Lecture Notes in Math. 1376 (1989), 105--112

\bibitem [T] {T}    M. Talagrand,
  {\em Embedding of $l_\infty^k$ and a theorem of Alon and Milman},
  GAFA (Israel, 1992--1994), Oper. Theory Adv. Appl., 77, Birkh\"{a}user, 
  Basel, 1995, 289--293. 

\bibitem [T-J] {T-J}    N. Tomczak-Jaegermann,
  {\em Banach-Mazur distances and finite dimensional operator ideals},
  Pitman, 1989

\bibitem [V] {V}    R. Vershynin,
  {\em Subsequences of frames},
  Submitted (1999)

\end{thebibliography}
}
\end{document}